\documentclass{amsart}
\usepackage{amssymb,amsmath,amsthm,amsbsy}
\usepackage[all]{xy}

\newtheorem{thm}{Theorem}[section]
\newtheorem{prop}[thm]{Proposition} 
\newtheorem{lemma}[thm]{Lemma}

\theoremstyle{definition} 
\newtheorem{dfn}[thm]{Definition}

\theoremstyle{remark}
\newtheorem{rmk}[thm]{Remark}
\newtheorem{ex}[thm]{Example}

\newtheorem*{ack}{Acknowledgements}

\newcommand{\la}[4]{
\xymatrix{#1 \ar[r] \ar@<2pt>[d] \ar@<-2pt>[d] & #2 \ar@<2pt>[d] \ar@<-2pt>[d] \\
#3 \ar[r] & #4}}

\newcommand{\arrows}{\rightrightarrows}
\newcommand{\defequal}{:=}

\newcommand{\reals}{{\mathbb R}}

\newcommand{\integers}{{\mathbb Z}}
\newcommand{\bitimes}[2]{\,_{#1}\!\!\times_{#2}}

\newcommand{\calO}{\mathcal{O}}

\newcommand{\vect}{\mathfrak{X}}

\renewcommand{\bigwedge}{\mbox{{\Large{$\wedge$}}}}
\renewcommand{\tilde}[1]{\widetilde{#1}}
\newcommand{\LA}{\mathcal{LA}}
\newcommand{\VB}{\mathcal{VB}}
\newcommand{\linv}[1]{\stackrel{\leftarrow}{#1}}
\newcommand{\rinv}[1]{\stackrel{\rightarrow}{#1}}
\newcommand{\tot}{\mathrm{tot}}
\newcommand{\id}{\mathrm{id}}
\newcommand{\tolabel}[1]{\stackrel{#1}{\to}}

\DeclareMathOperator{\ev}{ev}
\DeclareMathOperator{\rank}{rank}

\begin{document}
\title{$Q$-groupoids and their cohomology}
\author{Rajan Amit Mehta}
\address{Department of Mathematics\\
Washington University in Saint Louis\\
One Brookings Drive\\
Saint Louis, Missouri, USA 63130}
\email{raj@math.wustl.edu}

\begin{abstract}
We approach Mackenzie's $\LA$-groupoids from a supergeometric point of view by introducing \emph{$Q$-groupoids}, which are groupoid objects in the category of $Q$-manifolds. There is a faithful functor from the category of $\LA$-groupoids to the category of $Q$-groupoids.  We associate to every $Q$-groupoid a double complex that provides a model for the $Q$-cohomology of the classifying space.  As examples, we obtain models for equivariant $Q$- and orbifold $Q$-cohomology, and for equivariant Lie algebroid and orbifold Lie algebroid cohomology. We obtain double complexes associated to Poisson groupoids and groupoid-algebroid ``matched pairs''.
\end{abstract}

\maketitle

\section{Introduction}

 Mackenzie \cite{mac:dblie1, mac:dblie2} introduced $\LA$-groupoids as the intermediate objects between double Lie groupoids and double Lie algebroids.  An $\LA$-groupoid is essentially a groupoid object in the category of Lie algebroids.  Examples include the tangent prolongation groupoid
\begin{equation}\label{eqn:tangent}
\la{TG}{G}{TM}{M}
\end{equation}
 of a Lie groupoid $G \arrows M$ and the cotangent prolongation groupoid
\begin{equation}\label{eqn:cotangent}
\la{T^*G}{G}{A^*}{M}
\end{equation}
 of a Poisson groupoid $G \arrows M$ with Lie algebroid $A$.

The notion of $\LA$-groupoid may be considered a simultaneous generalization of those of Lie groupoids and Lie algebroids.  Indeed, if $G \arrows M$ is a Lie groupoid, then the square
\begin{equation}\label{trivgpd}
\la{G}{G}{M}{M},
\end{equation}
where the horizontal sides are trivial Lie algebroids, is an $\LA$-groupoid. Similarly, if $A \to M$ is a Lie algebroid, then the square
\begin{equation}\label{trivalgbd}
\la{A}{M}{A}{M},
\end{equation}
where the vertical sides are trivial Lie groupoids, is an $\LA$-groupoid.

In this paper, we approach $\LA$-groupoids from a supergeometric point of view in a way that extends Vaintrob's approach to Lie algebroids \cite{vaintrob}, in which the sheaf of differential algebras $(\bigwedge \Gamma(A^*), d_A)$ of a Lie algebroid $A \to M$ is interpreted as the function sheaf of a ($\integers$-)graded manifold equipped with a homological vector field.  This approach leads us to the notion of a \emph{$Q$-groupoid}\footnote{The terminology derives from the term \emph{$Q$-manifold}, due to Schwarz \cite{schwarz}, which refers to a supermanifold equipped with a homological vector field.}.  A $Q$-groupoid is a graded groupoid that is equipped with a compatible homological vector field.  There is a faithful functor, which we denote as $[-1]$, from the category of $\LA$-groupoids to the category of $Q$-groupoids.

It is known that a Lie groupoid $G$ gives rise to a simplicial manifold, the \emph{nerve} of $G$.  We show that a $Q$-groupoid gives rise to a simplicial $Q$-manifold.  From this simplicial $Q$-manifold, we obtain a natural double complex, which leads to a notion of $Q$-groupoid cohomology (and therefore of $\LA$-groupoid cohomology) that is a generalization of both Lie algebroid cohomology and Lie groupoid cohomology.  Of particular interest is the case of the tangent prolongation groupoid of a Lie groupoid $G$, whose double complex, known as the \emph{de Rham double complex} (see for example \cite{bss, dupont}) of $G$, is a model for the cohomology of the classifying space $BG$.  In general, $Q$-groupoid cohomology can be thought of as the $Q$-cohomology of the classifying space of the $Q$-groupoid, even though a classifying space may not exist in the category of $Q$-manifolds.

As an application, we introduce a notion of equivariant $Q$-cohomology, arising from a suitable action of a $Q$-group $\mathcal{G}$ on a $Q$-manifold $\mathcal{M}$.  In the case where the quotient $\mathcal{M}/\mathcal{G}$ is a graded manifold and the projection $\mathcal{M} \to \mathcal{M}/\mathcal{G}$ is a principal $\mathcal{G}$-bundle, we show that the equivariant $Q$-cohomology coincides with the $Q$-cohomology of the quotient.  As a special case, we give a noninfinitesimal model for equivariant Lie algebroid cohomology.

There is a close relationship between the supergeometric approach to $\LA$-group\-oids and Voronov's supergeometric approach to double Lie algebroids \cite{voronov:mack}.  This relationship will be detailed in \cite{m:qalg} and \cite{mw}.

The structure of the paper is as follows.  In \S\ref{sec:zgraded}, we give a brief introduction to $\integers$-graded manifolds.  In \S\ref{sec:q}, we define graded groupoids and $Q$-groupoids, and describe the $[-1]$ functor from the category of $\LA$-groupoids to the category of $Q$-groupoids.  In \S\ref{sec:cohomology}, we describe the functor from $Q$-groupoids to simplicial $Q$-manifolds and define $Q$-groupoid cohomology.

We then apply the $[-1]$ functor to some examples of $\LA$-groupoids and consider the resulting double complexes.  In \S\ref{sec:tangent}, we reproduce the de Rham double complex, which includes as special cases the simplicial model for equivariant cohomology and a model for orbifold cohomology.  In \S\ref{sec:equivq}--\ref{sec:orb}, we describe the more general notions of equivariant $Q$- and orbifold $Q$-cohomology.

\S\ref{sec:cotangent} addresses the cotangent prolongation groupoid of a Poisson groupoid, and in \S\ref{sec:vacant} we consider the case of vacant $\LA$-groupoids, which always arise from a matched pair consisting of a Lie groupoid and a Lie algebroid that compatibly act on each other.

\begin{ack}This work was partly supported by a fellowship from Conselho Nacional de Desenvolvimento Cient\'{i}fico e Tecnol\'{o}gico.  I would like to thank Alan Weinstein for providing much support and advice as my thesis advisor during the period when the main ideas of the paper were developed, as well as for comments on drafts of the paper.  I would also like to thank Henrique Bursztyn and Kirill Mackenzie for useful comments.  Finally, I would like to thank the anonymous referee, whose comments led to a substantially revised and greatly improved version of the paper.\end{ack}

\section{Geometry of graded manifolds}
\label{sec:zgraded}

Primarily to fix notation and terminology, we present here a brief introduction to $\integers$-graded manifolds.  For a more details introduction, see \cite{m:qalg} or \cite{mythesis}.  

Various notions of manifolds with $\integers$-gradings have appeared in the work of Kontsevich \cite{kont:deformation}, Roytenberg \cite{royt:graded}, \v{S}evera \cite{severa}, and Voronov \cite{voronov}, among others.  All of these notions share the common property of equipping a ($\integers_2$-graded) supermanifold with an additional $\integers$-grading.  

We prefer to give a definition that is intrinsic, in the sense that it does not rely on a prior definition of supermanifolds.  Since the $\integers_2$-grading in our definition is a consequence of the $\integers$-grading and not an antecedent structure, this definition is therefore less general than that of, e.g.\ \cite{voronov}, where independent $\integers$- and $\integers_2$-gradings are allowed\footnote{To include this case, one could easily extend the definition given here to $(\integers \times \integers_2)$-graded manifolds.  In fact, the even more general possibility of an $R$-grading, where $R$ is any ring equipped with a homomorphism to $\integers_2$, was already mentioned in the work of Kostant \cite{kostant}.}.  

Just as a manifold is locally modeled on $\reals^n$, a ($\integers$-)graded manifold is locally modeled on a \emph{graded coordinate space} $\reals^{\{p_i\}}$.

\subsection{Graded domains and graded manifolds}

Let $\{ p_i \}_{i \in \integers}$ be a nonnegative integer-valued sequence.  Denote by $\calO^{ \{p_i\} }$ the sheaf of graded, graded-commutative algebras on $\reals^{p_0}$ defined by
\begin{equation*}
\calO^{ \{p_i\} } (U) = C^\infty(U)\left[\bigcup_{i \neq 0} \{\xi_i^1, \dots, \xi_i^{p_i}\} \right]
\end{equation*}
for any open set $U \subseteq \reals^{p_0}$, where $\xi_i^k$ is of degree $-i$.  The coordinates $\{x^k, \xi_i^k\}$ are assumed to commute in the graded sense; in particular, $\xi_i^k \xi_j^\ell = (-1)^{ij}\xi_j^\ell\xi_i^k$.

\begin{dfn} The \emph{coordinate graded space} $\reals^{ \{p_i\} }$ is the pair $\left( \reals^{p_0} , \calO^{ \{p_i\} }\right)$.
\end{dfn}

\begin{rmk} The basic premise of graded geometry is that we treat $\reals^{ \{p_i\} }$ as if it were a space whose sheaf of ``smooth functions'' is $\calO^{ \{p_i\} }$.  Following this idea, we write $C^\infty\left(\reals^{ \{p_i\} }\right) \defequal \calO^{ \{p_i\} }$.
\end{rmk}

There is a natural surjection of sheaves $\ev: C^\infty\left(\reals^{ \{p_i\} }\right) \to C^\infty\left(\reals^{p_0}\right)$, called the \emph{evaluation map}, where the kernel is the ideal generated by all elements of nonzero degree.

\begin{dfn}
A \emph{graded domain} $\mathcal{U}$ of dimension $\{ p_i \}$ is a pair $\left(U, C^\infty(\mathcal{U})\right)$, where $U$ is an open subset of $\reals^{p_0}$ and $C^\infty(\mathcal{U}) \defequal \calO^{ \{p_i\} }|_U$.
A morphism of graded domains $\mu: \mathcal{U} \to \mathcal{V}$ consists of a smooth map $\mu_0: U \to V$ and a morphism of sheaves of graded algebras $\mu^*: C^\infty(\mathcal{V}) \to C^\infty(\mathcal{U})$ over $\mu_0$, such that $\ev \circ \mu^* = \mu_0^* \circ \ev$.
\end{dfn}

\begin{dfn}\label{dfn:supermfld}
A \emph{graded manifold} $\mathcal{M}$ of dimension $\{p_i\}$ is a pair $(M, C^\infty(\mathcal{M}))$, where $M$ (the \emph{support}) is a topological space and $C^\infty(\mathcal{M})$ is a sheaf of graded algebras (the \emph{sheaf of functions}) on $M$ that is locally isomorphic to a graded domain of dimension $\{p_i\}$.  A morphism of graded manifolds $\mu: \mathcal{M} \to \mathcal{N}$ consists of a map $\mu_0: M \to N$ and a sheaf morphism $\mu^*: C^\infty(\mathcal{N}) \to C^\infty(\mathcal{M})$ over $\mu_0$ that is locally a morphism of graded domains.
\end{dfn}

\begin{rmk}
It follows from Definition \ref{dfn:supermfld} that if $\mathcal{M} = (M, C^\infty(\mathcal{M}))$ is a dimension $\{p_i\}$ graded manifold, then the topological space $M$ automatically has the structure of a $p_0$-dimensional manifold.  The evaluation map describes an embedding of $M$ into $\mathcal{M}$.
\end{rmk}

\begin{rmk}
 Definition \ref{dfn:supermfld} allows for the possibility that $\mathcal{M}$ can be of infinite total dimension, in the sense that $\sum p_i = \infty$.  Nonetheless, a local coordinate description of any $f \in C^\infty(\mathcal{M})$ is, by definition, polynomial in the nonzero-degree coordinates.
\end{rmk}

In geometric situations, the graded manifolds of interest will often have function sheaves that are nonnegatively graded.
\begin{dfn}
A graded manifold $\mathcal{M}$ of dimension $\{p_i\}$ is said to be of \emph{degree d} if $p_i$ is nonzero only when $-d \leq i \leq 0$.
\end{dfn}

\subsection{The functor $[-1]$}

Let $E \to M$ be a vector bundle.  We may interpret the sheaf $\bigwedge \Gamma(E^*)$ as the function sheaf of a graded manifold, which we denote by $[-1]E$, with support $M$.  Its dimension is $\{p_i\}$, where $p_0 = \dim M$, $p_{-1} = \rank E$, and $p_i = 0$ for $i \neq 0,-1$.  In particular, $[-1]E$ is a degree $1$ graded manifold.

A bundle map
\begin{equation*}
\xymatrix{E \ar[d] \ar^\phi[r] & E' \ar[d] \\ M \ar[r] & M'}
\end{equation*}
induces a sheaf morphism $\phi^*: \bigwedge \Gamma(E'^*) \to \bigwedge \Gamma(E^*)$, which may be viewed as a morphism of graded manifolds $[-1]\phi: [-1]E \to [-1]E'$.  Thus we have a functor, denoted as $[-1]$, from the category of vector bundles to the category of (degree $1$) graded manifolds.

\begin{rmk}
Our notation differs from the notation in much of the literature, for example \cite{kont:deformation, royt:graded, severa, voronov}, where the graded manifold with function sheaf $\bigwedge \Gamma(E^*)$ is denoted as $E[1]$.  There are two separate distinctions at work here.  The first is that, in the spirit of supergeometry, we have interpreted the operation $[-1]$ to be a geometric operation as opposed to an algebraic one.  In other words, the \emph{fibres} of $[-1]E$ are of degree $-1$, whereas in the cited literature $E[1]$ is characterized by the property that the \emph{linear functions} are of degree $1$.  Because the degree of a homogeneous vector space is opposite in sign to the degree of its dual space, our degree shift operator differs by a sign from others' degree shift operator.  The second distinction, that the shift functor is placed on the left instead of the right, is of less importance in the present context of ordinary vector bundles; in the general setting of vector bundles in the category of graded manifolds, there exist both left and right shift functors, but their restrictions to purely even vector bundles are canonically isomorphic.
\end{rmk}

\subsection{Homological vector fields}
Let $\mathcal{M}$ be a graded manifold.
\begin{dfn} A \emph{vector field} of degree $j$ on $\mathcal{M}$ is a degree $j$ derivation $\phi$ of $C^\infty(\mathcal{M})$, i.e.\ a linear operator such that
\begin{equation*} |\phi f| = j + |f| \end{equation*}
and 
\begin{equation*} \phi(fg) = \phi(f)g + (-1)^{j|f|}f \phi(g) \end{equation*}
for any homogeneous functions $f, g \in C^\infty(\mathcal{M})$.
The space of degree $j$ vector fields on $M$ is denoted $\vect_j(\mathcal{M})$, and the space of all vector fields is $\vect(\mathcal{M}) \defequal \bigoplus_{j \in \integers} \vect_j(\mathcal{M})$.
\end{dfn}

The bracket $[\phi,\psi] \defequal \phi \psi - (-1)^{|\phi||\psi|}\psi\phi$ gives the space of vector fields on $\mathcal{M}$ the structure of a Lie superalgebra.  In particular, if $\phi$ is an odd degree vector field, then $[\phi,\phi] = 2\phi^2$ is not automatically zero.

\begin{dfn}A vector field $\phi$ on a graded manifold is called \emph{homological} if it is of odd degree and satisfies the equation $[\phi,\phi] = 0$.
\end{dfn}

\begin{dfn}[\cite{schwarz}] A \emph{$Q$-manifold} is a graded manifold equipped with a homological vector field.
\end{dfn}

\begin{rmk}
 For the examples of $Q$-manifolds in the present paper, the homological vector fields are of degree $1$.  When it is necessary to single out this special case, we will use the term \emph{$Q_1$-manifold} for a $Q$-manifold whose homological vector field is of degree $1$.
\end{rmk}

Let $(\mathcal{M}, \phi)$ and $(\mathcal{N}, \psi)$ be $Q$-manifolds.
\begin{dfn}
A $Q$-manifold morphism from $\mathcal{M}$ to $\mathcal{N}$ is a morphism of graded manifolds $\mu: \mathcal{M} \to \mathcal{N}$ such that $\phi$ and $\psi$ are $\mu$-related; that is, for all $f \in C^\infty(\mathcal{N})$, 
\begin{equation*}\mu^* (\psi f) = \phi (\mu^*f).\end{equation*}
\end{dfn}

\begin{ex}[\cite{vaintrob}]\label{onea}\label{algtoq}
Let $A \to M$ be a Lie algebroid.  Then the graded manifold $[-1]A$, equipped with the Lie algebroid differential $d_A$, is a $Q_1$-manifold.  Special cases include the odd tangent bundle $[-1]TM$ of a manifold $M$, equipped with the de Rham differential; the odd cotangent bundle $[-1]T^*M$ of a Poisson manifold $M$, equipped with the Lichnerowicz-Poisson differential; and $[-1]\mathfrak{g}$, where $\mathfrak{g}$ is a Lie algebra, equipped with the Chevalley-Eilenberg differential.  

On the other hand, any $Q_1$-manifold of degree $1$ is canonically of the form $[-1]A$ for some Lie algebroid $A$.  Lie algebroid morphisms from $A \to M$ to $A' \to M'$ are in one-to-one correspondence with $Q$-manifold morphisms from $[-1]A$ to $[-1]A'$, so the $[-1]$ functor gives an isomorphism from the category of Lie algebroids to the category of degree $1$ $Q_1$-manifolds.
\end{ex}

\section{Graded groupoids and $Q$-groupoids}
\label{sec:q}

\subsection{Fibre products}

Let $\mathcal{P}$ and $\mathcal{P}'$ be graded manifolds equipped with maps $\pi$ and $\pi'$, respectively, to a graded manifold $\mathcal{M}$.  The fibre product $\mathcal{P} \bitimes{\pi}{\pi'} \mathcal{P}'$ always exists in the category of ringed spaces, but is not necessarily a graded manifold.  As in the case of ordinary manifolds, a sufficient condition for $\mathcal{P} \bitimes{\pi}{\pi'} \mathcal{P}'$ to be a graded manifold is that $\pi$ be a submersion.

If $\mathcal{P}$, $\mathcal{P}'$, and $\mathcal{M}$ are $Q$-manifolds and the maps $\pi$ and $\pi'$ are $Q$-manifold morphisms, then the product homological vector field on $\mathcal{P} \times \mathcal{P'}$ preserves the ideal of functions that vanish on the fibre product $\mathcal{P} \bitimes{\pi}{\pi'} \mathcal{P}'$.  Therefore, when the fibre product is a graded manifold, it has a naturally-induced $Q$-manifold structure.  The fibre product $Q$-manifold structure is uniquely determined by the property that the natural projection maps to $\mathcal{P}$ and $\mathcal{P}'$ are $Q$-manifold morphisms.  If $\pi$ and $\pi'$ are surjective, then a converse statement holds.

\begin{lemma}\label{lemma:lift}  Suppose that $\pi$ is a surjective submersion and $\pi'$ is surjective, and let $p_1$ and $p_2$ be the natural projections of $\mathcal{P} \bitimes{\pi}{\pi'} \mathcal{P}'$ onto $\mathcal{P}$ and $\mathcal{P}'$, respectively.  Let $X$ and $X^\prime$ be vector fields on $\mathcal{P}$ and $\mathcal{P}^\prime$, respectively.  Then the following are equivalent:
\begin{enumerate}
\item There exists a vector field $Y$ on $\mathcal{M}$ that is both $\pi$-related to $X$ and $\pi^\prime$-related to $X^\prime$.
\item The product vector field $X \times X' \in \vect(\mathcal{P}\times \mathcal{P}')$ is tangent to the fibre product.
\item There is a (necessarily unique) vector field $\tilde{X}$ on $\mathcal{P} \bitimes{\pi}{\pi^\prime} \mathcal{P}^\prime$ that is $p_1$-related to $X$ and $p_2$-related to $X^\prime$.
\end{enumerate}
\end{lemma}
We omit the proof of Lemma \ref{lemma:lift}, as it retraces the case of ordinary manifolds.

We may apply the above discussion to Lie algebroids.  Let $A \to M$, $B \to P$, and $B' \to P'$ be Lie algebroids, and let $\phi: B \to A$ and $\phi': B' \to A$ be Lie algebroid morphisms over $\pi: P \to M$ and $\pi': P' \to M$, respectively.  Applying the $[-1]$ functor, we obtain $Q$-manifold morphisms $[-1]\phi$ and $[-1]\phi'$ from $[-1]B$ and $[-1]B'$, respectively, to $[-1]A$.  If $[-1]B \bitimes{[-1]\phi}{[-1]\phi'} [-1]B'$ is a graded manifold, then it has an induced $Q$-manifold structure such that the projections to $[-1]B$ and $[-1]B'$ are $Q$-manifold morphisms.  Translated back into the language of Lie algebroids, we have the following:

\begin{prop}\label{prop:fibprod1}
If the fibre product $B \bitimes{\phi}{\phi'} B' \to P \bitimes{\pi}{\pi'} P'$ is a vector bundle, then it has an induced Lie algebroid structure that is uniquely determined by the property that the projection maps to $B$ and $B'$ are Lie algebroid morphisms.
\end{prop}

\begin{rmk}
Higgins and Mackenzie \cite{higmac} stated Proposition \ref{prop:fibprod1} without proof, and Stefanini \cite{luca:la} recently gave a proof in ordinary (i.e.\ not supergeometric) language.
\end{rmk}

A sufficient condition for $[-1]B \bitimes{[-1]\phi}{[-1]\phi'} [-1]B'$ to be a graded manifold (and hence for $B \bitimes{\phi}{\phi'} B' \to P \bitimes{\pi}{\pi'} P'$ to be a vector bundle) is that $[-1]\phi$ be a submersion in the category of graded manifolds.  This is the case if and only if $\pi$ and the map $B \to P \times_M A$ are submersions; essentially, the degree $0$ coordinates on the fibres of $[-1]B$ over $[-1]A$ correspond to coordinates on the fibres of $P$ over $M$, while those of degree $1$ correspond to coordinates on the fibres of $B$ over $P \times_M A$.  Therefore, we have

\begin{prop}\label{prop:fibprod2}
If the maps $\pi$ and $B \to P \times_M A$ are submersions, then the fibre product $B \bitimes{\phi}{\phi'} B' \to P \bitimes{\pi}{\pi'} P'$ is a vector bundle and therefore inherits a fibre product Lie algebroid structure.
\end{prop}

\subsection{Graded groupoids}

A graded groupoid is a groupoid object in the category of graded manifolds.  In other words, a graded groupoid $\mathcal{G} \arrows \mathcal{M}$ is a pair of graded manifolds $(\mathcal{G}, \mathcal{M})$ equipped with surjective submersions $s,t : \mathcal{G} \to \mathcal{M}$ (source and target) and, letting $\mathcal{G}^{(2)}$ be the fibre product $\mathcal{G}\bitimes{s}{t}\mathcal{G}$, maps $m: \mathcal{G}^{(2)} \to \mathcal{G}$ (multiplication), $e: \mathcal{M} \to \mathcal{G}$ (identity), and $i: \mathcal{G}\to \mathcal{G}$ (inverse), satisfying a series of commutative diagrams that describe the various axioms of a groupoid.

\begin{rmk}If $\mathcal{G} \arrows \mathcal{M}$ is a graded groupoid, then there is an underlying ordinary groupoid $G \arrows M$, where $G$ and $M$ are the supports of $\mathcal{G}$ and $\mathcal{M}$, respectively.
\end{rmk}

We will now describe our main example of a graded groupoid.

\begin{ex}
 Consider a square
\begin{equation}\label{diag:vbg}
 \la{\Omega}{G}{E}{M},
\end{equation}
where the horizontal sides are vector bundles and the vertical sides are groupoids.  Let $p_\Omega: \Omega \to G$ and $p_E: E \to M$ be the projection maps, let $s$, $t$ and $m$ be the source, target, and multiplication maps for $G \arrows M$, and let $\tilde{s}$, $\tilde{t}$, and $\tilde{m}$ be the source, target, and multiplication maps for $\Omega \arrows E$.  If $(\tilde{s} \times p_\Omega) : \Omega \to E \bitimes{p_E}{s} G$ is a submersion, then $\Omega^{(2)} \to G^{(2)}$ is a vector bundle; see Proposition \ref{prop:fibprod2}.  If, furthermore, the groupoid structure maps for $\Omega \arrows E$  all form vector bundle morphisms over the corresponding maps for $G \arrows M$, then (\ref{diag:vbg}) is called a \emph{$\VB$-groupoid} \cite{mac:dblie1, pradines2}.  In this case, the $[-1]$ functor may be applied to the vector bundles to obtain the graded groupoid $[-1]\Omega \arrows [-1]E$.
\end{ex}

\subsection{Multiplicative vector fields}

Let $\mathcal{G} \arrows \mathcal{M}$ be a graded groupoid.  Let $p_1$ and $p_2$ be the natural maps from $\mathcal{G}^{(2)}$ to $\mathcal{G}$ that project onto the first and second components, respectively.  We may use the three maps $m, p_1, p_2 : \mathcal{G} \to \mathcal{G}^{(2)}$ to give a ``simplical''\footnote{The three maps $p_2$, $m$, and $p_1$ are face maps for the simplicial (graded) manifold associated to $\mathcal{G}$.  This is discussed in more detail in \S\ref{sec:cohomology}.} definition of multiplicative vector fields:

\begin{dfn}\label{dfn:mult}
A vector field $\psi$ on $\mathcal{G}$ is \emph{multiplicative} if there exists a vector field $\psi^{(2)}$ on $\mathcal{G}^{(2)}$ that is $p_1$-, $p_2$-, and $m$-related to $\psi$.
\end{dfn}

\begin{rmk}\label{rmk:mult}
By Lemma \ref{lemma:lift}, a multiplicative vector field could be equivalently defined as a pair of vector fields $\psi \in \vect(\mathcal{G})$, $\psi^{(0)} \in \vect(\mathcal{M})$, such that $\psi$ is $s$- and $t$-related to $\psi^{(0)}$, and such that the restriction of the product vector field $\psi \times \psi$ to $\mathcal{G}^{(2)}$ is $m$-related to $\psi$.  The induced vector field $\psi^{(0)}$ is called the \emph{base vector field} of $\psi$.
\end{rmk}

\begin{lemma}\label{lemma:multbase}
If $\psi$ is a multiplicative vector field with base vector field $\psi^{(0)}$, then $\psi$ is $e$-related to $\psi^{(0)}$.
\end{lemma}
\begin{proof}
Let $\Delta_0^1$ and $\Delta_1^1$ be the degeneracy maps from $\mathcal{G}$ to $\mathcal{G}^{(2)}$, defined componentwise by $\Delta_0^1(g) = ((e \circ t)(g),g)$ and $\Delta_1^1(g) = (g, (e \circ s)(g))$.  It follows from the groupoid axioms that $m \circ \Delta_i^1 = \id$ for $i = 0,1$.

If $\psi$ is a multiplicative vector field, then 
\begin{equation}\label{eqn:multbase1}
(\Delta_1^1)^* \circ \psi^{(2)} \circ m^* = (\Delta_1^1)^* \circ m^* \circ \psi = \psi = \psi \circ (\Delta_1^1)^* \circ m^*.
\end{equation}
Similarly,
\begin{equation}\label{eqn:multbase2}
(\Delta_1^1)^* \circ \psi^{(2)} \circ p_1^* = \psi \circ (\Delta_1^1)^* \circ p_1^*.
\end{equation}

The axioms for the inverse map imply that the two maps
\begin{align*}
 \begin{split}
  \mathcal{G}^{(2)} &\to \mathcal{G} \bitimes{t}{t}\mathcal{G}, \\
  (g,h) &\mapsto (g,gh), \\
 \end{split}
\begin{split}
   \mathcal{G}^{(2)} &\to \mathcal{G} \bitimes{s}{s}\mathcal{G}, \\
  (g,h) &\mapsto (gh,h), \\
\end{split}
\end{align*}
are invertible.  Therefore, by the uniqueness property of Lemma \ref{lemma:lift}, equations \eqref{eqn:multbase1} and \eqref{eqn:multbase2} imply that $\psi$ and $\psi^{(2)}$ are $\Delta_1^1$-related, and it follows that
\begin{equation}\label{eqn:multbase3}
 (\Delta_1^1)^* \circ \psi^{(2)} \circ p_2^* = \psi \circ (\Delta_1^1)^* \circ p_2^*.
\end{equation}
The left hand side of \eqref{eqn:multbase3} is $(\Delta_1^1)^* \circ p_2^* \circ  \psi  = s^* \circ e^* \circ \psi$, whereas the right hand side is $\psi \circ s^* \circ e^* = s^* \circ \psi^{(0)} \circ e^*$.  Since $s^*$ is injective, we conclude that $e^* \circ \psi = \psi^{(0)} \circ e^*$.
\end{proof}

\subsection{$Q$-groupoids}

Let us recall the definition of an $\LA$-groupoid.  Consider a square
\begin{equation}\label{eq:la2}
\la{\Omega}{G}{A}{M},
\end{equation}
where the horizontal sides are Lie algebroids and the vertical sides are Lie groupoids.  Let $p_\Omega: \Omega \to G$ and $p_A: A \to M$ be the projection maps, let $s$, $t$ and $m$ be the source, target, and multiplication maps for $G \arrows M$, and let $\tilde{s}$, $\tilde{t}$, and $\tilde{m}$ be the source, target, and multiplication maps for $\Omega \arrows A$.

\begin{dfn}\label{dfn:la}
The square (\ref{eq:la2}) is an $\LA$-groupoid if

\begin{enumerate}
\item The maps $\tilde{s}$ and $\tilde{t}$ from $\Omega$ to $A$ form Lie algebroid morphisms over the corresponding maps from $G$ to $M$.
\item The map $(\tilde{s} \times {p_\Omega}) : \Omega \to A \bitimes{p_A}{s} G$ is a surjective submersion.
\item The map $\tilde{m}: \Omega^{(2)} \to \Omega$ is a Lie algebroid morphism over $m: G^{(2)} \to G$.
\end{enumerate}
\end{dfn}

\begin{rmk}
In the third condition of Definition \ref{dfn:la}, $\Omega^{(2)} \to G^{(2)}$ is equipped with the fibre product Lie algebroid structure, which exists as a consequence of the first two statements; see Proposition \ref{prop:fibprod2}.
\end{rmk}

Translating Definition \ref{dfn:la} into the language of $Q$-manifolds, we immediately obtain the following:

\begin{thm} \label{thm:latoq}
A $\VB$-groupoid of the form (\ref{eq:la2}), where the horizontal sides are Lie algebroids, is an $\LA$-groupoid if and only if the associated homological vector field on the graded groupoid $[-1]\Omega \arrows [-1]A$ is multiplicative.
\end{thm}

\begin{dfn} A \emph{$Q$-groupoid} is a graded groupoid equipped with a multiplicative, homological vector field.
\end{dfn}

\begin{rmk}The $[-1]$ functor from the category of $\LA$-groupoids to the category of $Q$-groupoids is faithful and full, so the statement of Example \ref{algtoq} directly generalizes as follows:  \emph{The $[-1]$ functor is an isomorphism of categories from the category of $\LA$-groupoids to the category of degree $1$ $Q_1$-groupoids.}
\end{rmk}

By applying Theorem \ref{thm:latoq} to the tangent prolongation groupoid (\ref{eqn:tangent}) and the cotangent prolongation groupoid (\ref{eqn:cotangent}), we immediately obtain the examples $[-1]TG \arrows [-1]TM$, when $G \arrows M$ is a groupoid, and $[-1]T^*G \arrows [-1]A^*$, when $G \arrows M$ is a Poisson groupoid.  These examples are described in more detail in \S\ref{sec:tangent} and \S\ref{sec:cotangent}, respectively.

Interesting examples of $Q$-groupoids that are not of degree $1$ may arise from Courant algebroids.  The correspondence between Courant algebroids and degree $2$ symplectic $Q_1$-manifolds, due to Roytenberg \cite{royt:graded} and \v{S}evera \cite{severa}, leads us to consider the following example of a degree $2$ symplectic $Q$-groupoid.

\begin{ex}
\label{ex:courant}
Let $G \arrows M$ be a Lie groupoid with Lie algebroid $A$.  Then $[-1]T^*G \arrows [-1]A^*$ is a degree $1$ symplectic groupoid.  
We may apply the $[-1]T$ functor to obtain the $Q$-groupoid $[-1]T([-1]T^*G) \arrows [-1]T([-1]A^*)$.  There is a degree $2$ symplectic structure on $[-1]T([-1]T^*G)$, arising from the canonical symplectic structure on $T^*G$, and, similarly, there is a Poisson structure on $[-1]T([-1]A^*)$ arising from the linear Poisson structure on $A^*$.  The relation to Courant algebroids is that $[-1]T([-1]T^*G)$ is the degree $2$ symplectic $Q$-manifold corresponding to the standard Courant algebroid $TG \oplus T^*G$.  Although $[-1]T([-1]A^*)$ is only Poisson, we may associate it to a structure on $A \oplus T^*M$ similar to that of a Courant algebroid, except that the bilinear form may be degenerate.  We observe that $TG \oplus T^*G$ has a natural Lie groupoid structure over the dual bundle $TM \oplus A^*$, and we suggest that this is a first example of a more general notion that might be called ``Courant groupoids''.
\end{ex}

\section{Cohomology of $Q$-groupoids}
\label{sec:cohomology}

Let $G \arrows M$ be a Lie groupoid.  It is well known that there is an associated simplicial manifold
\begin{equation*}
\xymatrix{
\cdots \ar@<6pt>[r]\ar@<2pt>[r]\ar@<-2pt>[r]\ar@<-6pt>[r] & G^{(2)} \ar@<4pt>[r]\ar[r]\ar@<-4pt>[r] & G^{(1)} = G \ar@<2pt>[r]\ar@<-2pt>[r] & G^{(0)} = M,
}
\end{equation*}
known as the \emph{nerve} $B_\bullet G$ of $G$, where $G^{(q)}$ is the manifold of compatible $q$-tuplets of elements of $G$.  If $\mathcal{G} \arrows \mathcal{M}$ is a graded Lie groupoid, then the nerve construction may still be carried out, yielding a simplicial graded manifold; we shall now describe this construction more explicitly.

Let $\mathcal{G}^{(0)} = \mathcal{M}$ and 
\begin{equation*}
\mathcal{G}^{(q)} = \underbrace{ \mathcal{G} \bitimes{s}{t} \cdots \bitimes{s}{t} \mathcal{G}}_{q}
\end{equation*}
for $q > 0$.  Here, the fibre products exist as graded manifolds since $s$ and $t$ are assumed to be submersions.  The face maps $\sigma_i^q : \mathcal{G}^{(q)} \to \mathcal{G}^{(q-1)}$ are given by $\sigma_0^1 = s$, $\sigma_1^1 = t$, and
\begin{align*}
\sigma_0^q &= p_2 \times \id \times \cdots \times \id, \\
\sigma_i^q &= \underbrace{\id \times \cdots \times \id}_{i - 1} \times m \times \underbrace{\id \times \cdots \id}_{q - i - 1}, & 0 < i < q, \\
\sigma_q^q &= \id \times \cdots \times \id \times p_1
\end{align*}
for $q > 1$.  The degeneracy maps are given by $\Delta_0^0 = e$ and 
\begin{align*}
\Delta_i^q &= \underbrace{\id \times \cdots \times \id}_{i} \times \Delta_0^1 \times \underbrace{\id \times \cdots \id}_{q - i - 1}, & i < q, \\
\Delta_i^q &= \underbrace{\id \times \cdots \times \id}_{i - 1} \times \Delta_1^1 \times \underbrace{\id \times \cdots \id}_{q - i}, & 0 < i
\end{align*}
for $q > 0$.  One may verify that the simplicial relations hold as a consequence of the groupoid axioms.

Now suppose that $\mathcal{G}$ is a $Q$-groupoid with homological vector field $\psi$.  Since $\psi$ is multiplicative, there exists by Remark \ref{rmk:mult} a base vector field $\psi^{(0)}$ on $\mathcal{G}^{(0)}$ that is $\sigma_0^1$- and $\sigma_1^1$-related to $\psi$.  By Definition \ref{dfn:mult}, there exists a lift $\psi^{(2)}$ on $\mathcal{G}^{(2)}$ that is $\sigma_i^2$-related to $\psi$ for $i = 0,1,2$. Furthermore, since the higher face maps $\sigma_i^q$ for $q>2$, are defined in terms of $\sigma_i^2$, there exist higher lifts $\psi^{(q)}$ such that $\psi^{(q)}$ and $\psi^{(q-1)}$ are $\sigma_i^q$-related for all $0 \leq i \leq q$.  As a consequence of Lemma \ref{lemma:multbase}, it is also the case that $\psi^{(q)}$ and $\psi^{(q+1)}$ are $\Delta_i^q$-related for all $0 \leq i \leq q$.  In summary:

\begin{thm}\label{thm:qgpd}
  If $\mathcal{G}$ is a $Q$-groupoid, then the nerve $B_\bullet \mathcal{G}$ is a simplicial $Q$-manifold.
\end{thm}

An immediate consequence of Theorem \ref{thm:qgpd} (specifically, the fact that $\psi^{(q)}$ and $\psi^{(q-1)}$ are $\sigma_i^q$-related) is that the action of the vector fields $\psi^{(q)}$ as derivations commutes with the groupoid coboundary operator $\delta^q: C^\infty(\mathcal{G}^{(q)}) \to C^\infty(\mathcal{G}^{(q+1)})$, defined as
\begin{equation*}
\delta^q = \sum_{i=0}^{q+1} (-1)^i \left(\sigma_i^{q+1}\right)^*.
\end{equation*}

Let $C^{p,q}(\mathcal{G}) \defequal C^\infty_p(\mathcal{G}^{(q)})$ denote the space of degree $p$ functions on $\mathcal{G}^{(q)}$.  

\begin{dfn}
\begin{enumerate}
\item The \emph{standard cochain complex} of the $Q$-groupoid $(\mathcal{G}, \psi)$
is the double complex $\left(C^{p,q}(\mathcal{G}), \delta, \psi\right)$.
\item The \emph{$Q$-groupoid cohomology} $H_\psi^\bullet(\mathcal{G})$ is the total cohomology  $H^\bullet \left(C_{\tot}^\bullet(\mathcal{G}), D\right)$, where $C_\tot^r(\mathcal{G}) \defequal \bigoplus_{p+q = r} C^{p,q}(\mathcal{G})$ and $D \defequal \psi + (-1)^p \delta$.
\end{enumerate}
\end{dfn}

In the case where $\mathcal{G}$ arises from an $\LA$-groupoid
\begin{equation*}
\la{\Omega}{G}{A}{M},
\end{equation*}
then we may refer to $H_\psi(\mathcal{G})$ as the \emph{$\LA$-groupoid cohomology} of $\Omega$.  Let us first consider the trivial examples.  

\begin{ex}\label{ex:trivgpd}
If $G \arrows M$ is a Lie groupoid, then the $Q$-groupoid that arises from the $\LA$-groupoid (\ref{trivgpd}) is simply $G \arrows M$ with the zero homological vector field.  In this case, $C^{p,q}(G) = 0$ for $p \neq 0$, and the total complex may be directly identified with the smooth Eilenberg-Maclane complex of $G$.
\end{ex}

\begin{ex}\label{ex:trivalgbd}
If $A \to M$ is a Lie algebroid, then the $Q$-groupoid arising from the $\LA$-groupoid (\ref{trivalgbd}) is $[-1]A \arrows [-1]A$ with the homological vector field $d_A$.  In this case, $C^{p,q}([-1]A) = \bigwedge^p \Gamma (A^*)$ for all $q$.  Since $\sigma_i^q = \id$ for all $i$, we have that $\delta^q = 0$ for even $q$ and $\delta^q = \id$ for odd $q$.  At the first stage, the spectral sequence for this double complex collapses to the Lie algebroid cohomology complex $\left(\bigwedge \Gamma(A^*), d_A\right)$.
\end{ex}

From Examples \ref{ex:trivgpd} and \ref{ex:trivalgbd}, we see that $\LA$-groupoid cohomology generalizes both Lie algebroid cohomology and Lie groupoid cohomology.  In the sections that follow, we will consider more interesting examples.

\section{Examples of $Q$-groupoid cohomology}

\subsection{The de Rham double complex}\label{sec:tangent}

Let $G \arrows M$ be a Lie groupoid.  We apply the $[-1]$ functor to the tangent prolongation groupoid (\ref{eqn:tangent}) and obtain the $Q$-groupoid $[-1]TG \arrows [-1]TM$, where the homological vector field is the de Rham differential $d$.

For each $q$, there is a natural identification of $([-1]TG)^{(q)}$ with $[-1]T(G^{(q)})$.  Thus the space of cochains is
\begin{equation*}
	C^{p,q}([-1]TG) = \Omega^p(G^{(q)}).
\end{equation*}
The double complex $(\Omega^p(G^{(q)}), \delta, d)$ is known as the \emph{de Rham double complex} \cite{ltx} of $G$.  The de Rham double complex is a special case of the de Rham complex of a simplicial manifold \cite{bss}, which is a model for the cohomology of the geometric realization.  In this case, the simplicial manifold is the nerve of $G$, whose geometric realization is the classifying space $BG$; see \cite{segal}.  Therefore, the $\LA$-groupoid cohomology of $TG$ is equal to $H^\bullet (BG; \reals)$.

Since the de Rham double complex is already well-known, $\LA$-groupoid cohomology does not provide any new information about classifying spaces.  However, as we illustrate in \S\ref{sec:equivariant} and \S\ref{sec:orb}, the $\LA$-groupoid point of view may be used to produce interesting generalizations of the de Rham double complex.

\begin{ex}\label{ex:equivariant}
Let $\Gamma$ be a Lie group that acts (from the right) on a manifold $M$.  Then it may be shown that the classifying space of the action groupoid $M \times \Gamma \arrows M$ is the homotopy quotient $M \times E\Gamma/\Gamma$.  The de Rham double complex of the action groupoid therefore computes the equivariant cohomology $H^\bullet_\Gamma(M)$.
\end{ex}

\begin{ex}\label{ex:orb}Let $G \arrows M$ be an groupoid presenting an orbifold $X$.  It is known (see, e.g.\ \cite{mp:orb2}) that the orbifold cohomology is isomorphic to the cohomology of $BG$.  Therefore, the de Rham double complex of $G$ is a model for the orbifold cohomology of $X$.
\end{ex}

\begin{ex}
The double complex of $[-1]T([-1]T^*G) \arrows [-1]T([-1]A^*)$, the $Q$-groupoid from Example \ref{ex:courant}, is the de Rham double complex of the graded groupoid $[-1]T^*G \arrows [-1]A^*$.  Since this graded groupoid has a linear structure over the ordinary groupoid $G \arrows M$, the double complex retracts\footnote{In general for a degree $1$ graded groupoid $\mathcal{G} \arrows \mathcal{M}$, there exists a natural linear structure over the underlying ordinary groupoid $G \arrows M$, and it follows that $H^\bullet(B\mathcal{G}) = H^\bullet(BG)$.  However, it seems plausible that one could construct a graded groupoid (possibly with both positively- and negatively-graded coordinates) for which there does not exist a linear structure over the underlying ordinary groupoid, and for which $H^\bullet(B\mathcal{G})$ does not equal $H^\bullet(BG)$.} to the de Rham double complex of $G$, and the $Q$-groupoid cohomology is just $H^\bullet(BG)$.  It may be interesting to see if one can ``twist'' the homological vector field on $[-1]T([-1]T^*G)$ by introducing a closed $3$-form on $G$.
\end{ex}

\subsection{Equivariant $Q$-cohomology}
\label{sec:equivq}

In Example \ref{ex:equivariant}, we learned that equivariant de Rham cohomology forms a special case of $Q$-groupoid cohomology.  In this section, we will extend this idea by introducing the notion of equivariant $Q$-cohomology.

Let $(\mathcal{M}, \phi)$ be a $Q$-manifold, and let $(\mathcal{G}, \psi)$ be a $Q$-group (that is, a $Q$-groupoid over a point) with multiplication map $\mu: \mathcal{G} \times \mathcal{G} \to \mathcal{G}$.

\begin{dfn}
 A (right) action of $\mathcal{G}$ on $\mathcal{M}$ is a \emph{$Q$-action} if the action map $\mathcal{M} \times \mathcal{G} \to \mathcal{M}$ is a morphism of $Q$-manifolds.
\end{dfn}

Given a right action of $\mathcal{G}$ on $\mathcal{M}$, we can form the action groupoid $\mathcal{M} \times \mathcal{G} \arrows \mathcal{M}$ in the same manner as in the ordinary case; specifically, the source map $s$ is the action map, the target map $t$ is projection onto $\mathcal{M}$, and, under the identification of $(\mathcal{M} \times \mathcal{G})^{(2)}$ with $\mathcal{M} \times \mathcal{G} \times \mathcal{G}$, the multiplication map is $m \defequal \id \times \mu$.  Clearly, $t$ and $\mu$ are automatically $Q$-manifold morphisms, so the following statement is immediate.

\begin{prop}\label{prop:qaction}
 The action groupoid $\mathcal{M} \times \mathcal{G} \arrows \mathcal{M}$ is a $Q$-groupoid if and only if the action of $\mathcal{G}$ on $\mathcal{M}$ is a $Q$-action.
\end{prop}

Generalizing the notion of equivariant de Rham cohomology (c.f.\ Example \ref{ex:equivariant}), we make the following definition.

\begin{dfn}
 Let $\mathcal{M}$ be a $Q$-manifold, and let $\mathcal{G}$ be a $Q$-group with a $Q$-action on $\mathcal{M}$.  The \emph{equivariant $Q$-cohomology} of $\mathcal{M}$ is the $Q$-groupoid cohomology of the action groupoid.
\end{dfn}
In analogy with the ordinary case, we may define the $Q$-cohomology of $B\mathcal{G}$ as the equivariant $Q$-cohomology of a point.  Again taking a clue from the ordinary case, we should expect the equivariant $Q$-cohomology to be the $Q$-cohomology of the quotient if the action is free and proper.  We will prove such a result, but first we need to define free and proper actions in the category of graded manifolds.

The algebra of $\mathcal{G}$-invariant functions on $\mathcal{M}$ is 
\begin{equation*}
 C^\infty_\mathcal{G}(\mathcal{M}) \defequal \{ f \in C^\infty(\mathcal{M}) : s^*(f) = t^*(f)\}.
\end{equation*}
Let $M$ and $G$ be the supports of $\mathcal{M}$ and $\mathcal{G}$, respectively.  Then $\left(M/G, C^\infty_\mathcal{G}(\mathcal{M})\right)$ is the ringed space corresponding to the quotient $\mathcal{M}/\mathcal{G}$.  Of course, the quotient may not be a graded manifold.  Still, if $\mathcal{M}/\mathcal{G}$ is a graded manifold, then it inherits the structure of a $Q$-manifold, since the homological vector field $\phi$ preserves the subalgebra of invariant functions.

We will say that the action of $\mathcal{G}$ on $\mathcal{M}$ is free and proper if $\mathcal{M}/\mathcal{G}$ is a graded manifold and the projection $\mathcal{M} \to \mathcal{M}/\mathcal{G}$ is a principal $\mathcal{G}$-bundle.

\begin{thm}\label{thm:equivfree}
 Suppose $\mathcal{G}$ is a $Q$-group with a free and proper $Q$-action on a $Q$-manifold $\mathcal{M}$.  Then the equivariant $Q$-cohomology of $\mathcal{M}$ is equal to the $Q$-cohomology of $\mathcal{M}/\mathcal{G}$.
\end{thm}

\begin{proof}
We will compute the groupoid cohomology of the action groupoid, and we will see that it is equal to $C^\infty(\mathcal{M}/\mathcal{G})$.  We will then be able to conclude that the spectral sequence associated to the standard complex of the action groupoid collapses at the first stage to the $Q$-complex of $\mathcal{M}/\mathcal{G}$.

First, let us consider the case where $\mathcal{M} = \mathcal{G}$, and where the action of $\mathcal{G}$ on itself is by right-multiplication.  In this case, the space of $(p,q)$-cochains in the standard complex is $C_p^\infty(\mathcal{G}^{q+1})$.  Let $\tau^q: 
C^\infty(\mathcal{G}^{q+2}) \to C^\infty(\mathcal{G}^{q+1})$ be given by evaluation at the identity element in the first component, and let $\delta^q: C^\infty(\mathcal{G}^{q+1}) \to C^\infty(\mathcal{G}^{q+2})$ be the groupoid coboundary operator.  If $q>0$, then $\tau^q \delta^q = \id - \delta^{q-1} \tau^{q-1}$, and it follows that the groupoid cohomology is trivial for $q>0$.  Additionally, $\tau^0 \delta^0 = \id - \ev_e$, where $\ev_e(f)$ is the constant function whose value is $f(e)$ for all $f \in C^\infty(\mathcal{G})$.  Thus, the $\delta^0$-closed functions on $\mathcal{G}$ are the constant functions, so the groupoid cohomology in degree $0$ is $\reals$.

Second, consider the case where $\mathcal{M} = \mathcal{U} \times \mathcal{G}$, where $\mathcal{U}$ is any graded manifold, and where $\mathcal{G}$ acts on the right component by right-multiplication.  Then, extending the previous case, we let $\tilde{\tau}^q = \id \times \tau^q: C^\infty(\mathcal{U} \times \mathcal{G}^{q+2}) \to C^\infty(\mathcal{U} \times \mathcal{G}^{q+1})$.  The groupoid coboundary operator in this case is $\tilde{\delta} = \id \times \delta$, so by the same argument as in the previous case, we have that the groupoid cohomology is trivial for $q>0$ and equal to $C^\infty(\mathcal{U})$ in degree $0$.

Finally, we turn to the general case.  Consider the projection
\begin{equation}\label{eqn:quotient}
\mathcal{M} \times \mathcal{G}^q \tolabel{\pi_q} \mathcal{M}/\mathcal{G}.
\end{equation}
Under the assumption that the action is free and proper, there exists an open cover $\{U_\alpha\}$ of the support $M/G$  of $\mathcal{M}/\mathcal{G}$, such that the restriction of \eqref{eqn:quotient} to each $U_\alpha$ is isomorphic (as a simplicial graded manifold) to $\mathcal{U}_\alpha \times \mathcal{G}^{q+1} \to \mathcal{U}_\alpha$, where $\mathcal{U}_\alpha$ is the restriction of $\mathcal{M}/\mathcal{G}$ to $U_\alpha$.  Thus, as we have seen in the previous case, the groupoid cohomology of the action groupoid locally retracts to $C^\infty(\mathcal{M}/\mathcal{G})$, and it remains to show that the local retracts can be patched together.

Let $\{\varphi_\alpha\}$ be a partition of unity subordinate to $\{U_\alpha\}$.  This partition of unity may be pulled back to partitions of unity $\{\pi_q^*(\varphi_\alpha)\}$ on $M \times G^q$, subordinate to $\{\pi_q^{-1}(U_\alpha)\}$.  Furthermore, these pulled back partitions of unity are preserved by the face maps.  A direct computation then shows that such a partition of unity can be used to patch together the local retracts at the level of groupoid cohomology.  Therefore, the groupoid cohomology of the action groupoid equals $C^\infty(\mathcal{M}/\mathcal{G})$, and by the spectral sequence argument given at the beginning of the proof, we conclude that the $Q$-groupoid cohomology of the action groupoid equals the $Q$-cohomology of the quotient.
\end{proof}

\subsection{Equivariant Lie algebroid cohomology}\label{sec:equivariant}

This cohomology was originally introduced by Ginzburg \cite{ginzburg} as a natural generalization of his theory of \emph{equivariant Poisson cohomology}.  More recently, Bruzzo, et al.\ \cite{bruzzo} have proven a corresponding localization theorem.  The model introduced by Ginzburg is a generalization of the Cartan model, which is in terms of the infinitesimal data of the action.  In this section, we introduce, as a special case of equivariant $Q$-cohomology, a model for equivariant Lie algebroid cohomology that is noninfinitesimal and is therefore useful, for example, in the case of a discrete group action.

Let $A \to M$ be a Lie algebroid and let $G$ be a Lie group.  We equip $A \times TG \to M \times G$ with the product Lie algebroid structure.

\begin{dfn}\label{dfn:gpaction}
An \emph{$A$-action} of $G$ is a (right) action of $TG$ on $A$ such that the action map $\widetilde{s}: A \times TG \to A$ is an algebroid morphism.
\end{dfn}

\begin{rmk}
An \emph{$A$-action} possesses an underlying action map $s: M \times G \to M$.  In practice, one will often begin with an action of $G$ on $M$ which one will then seek to extend to an $A$-action.
\end{rmk}

From the perspective of graded geometry, an $A$-action of $G$ is a $Q$-action of $[-1]TG$ on $[-1]A$.  

\begin{dfn}
Let $\widetilde{s}: A \times TG \to A$ be an $A$-action.  The \emph{equivariant Lie algebroid cohomology} of $A$ is the equivariant $Q$-cohomology of $[-1]A$ with respect to the action of $[-1]TG$.
\end{dfn}

\begin{ex}
For $s: M \times G \to M$ an action map, $Ts: TM \times TG \to TM$ is the unique $TM$-action that lifts $s$.  In this case, the equivariant Lie algebroid cohomology is equal to the equivariant cohomology $H^\bullet_G(M)$.
\end{ex}

\begin{ex}
If $G$ is a discrete group, then any action $\widetilde{s}: A \times G \to A$ where $G$ acts by Lie algebroid automorphisms is an $A$-action.  The resulting double complex is then of the form $\left(\bigwedge^q \Gamma(A^*) \times G^q, d_A, \delta\right)$, where $\delta$ is the coboundary operator for group cohomology with coefficients in $\bigwedge^q \Gamma(A^*)$.  In the case where $G$ is finite, the fact that $H^n (G;\bigwedge^q \Gamma(A^*))$ vanishes for $n>0$ implies that the spectral sequence collapses to the complex $\left(\left(\bigwedge^q \Gamma(A^*)\right)^G, d_A\right)$ of invariant Lie algebroid forms.
\end{ex}

The following theorem is a special case of Theorem \ref{thm:equivfree}.
\begin{thm}
 Let $A \to M$ be a Lie algebroid, and let $G$ be a Lie group with an $A$-action such that the action of $TG$ on $A$ is free and proper.  Then the quotient $A/TG \to M/G$ naturally inherits a Lie algebroid structure, and the equivariant Lie algebroid cohomology of $A$ equals the Lie algebroid cohomology of $A/TG$.
\end{thm}

\subsection{$Q$-orbifolds and Lie algebroid structures over orbifolds}\label{sec:orb}

Example \ref{ex:orb} gives a motivation to use the $Q$-groupoid point of view to define the notion of \emph{$Q$-orbifold}.  Specifically, if $G \arrows M$ is an orbifold groupoid presenting an orbifold $X$, then a $Q$-orbifold over $X$ is presented\footnote{Strictly speaking, one should check that this notion is well-defined, in the sense of respecting Morita equivalences (see e.g.\ \cite{moerdijk:orbgpds}).  We have not made an attempt to show this, and as such, the ideas in this section should be considered tentative.} by a $Q$-groupoid $\mathcal{G} \arrows \mathcal{M}$ whose support is $G \arrows M$.  As a special case, an $\LA$-groupoid
\begin{equation*}
\la{\Omega}{G}{A}{M}
\end{equation*}
may be viewed as representing a Lie algebroid over $X$.  In these cases, the $Q$-groupoid cohomology may be interpreted as ``orbifold $Q$-cohomology'' or ``orbifold Lie algebroid cohomology''.

\begin{ex}
One possible way to describe a Poisson structure on an orbifold $X$ would be to give an \'{e}tale groupoid representative $G \arrows M$ for which $G$ and $M$ are Poisson manifolds and the structure maps are Poisson morphisms.  Orbifold Poisson cohomology can then be defined as the cohomology of the $\LA$-groupoid\footnote{Note that when $G \arrows M$ is an \'{e}tale groupoid, there naturally exists a cotangent groupoid $T^*G \arrows T^*M$.  This is clearly different from the cotangent prolongation groupoid (\ref{eqn:cotangent}), which exists even when $G \arrows M$ is not \'{e}tale.}
\begin{equation*}
\la{T^*G}{G}{T^*M}{M}.
\end{equation*}	
In this case, the standard cochain complex is $(\vect^p(G^{(q)}), \delta, d_\pi)$.
\end{ex}

\begin{ex}If $X = M/\Gamma$ is the global quotient of a manifold $M$ by the action of a finite group $\Gamma$, then $X$ may be represented by the action groupoid $M \times \Gamma \arrows M$.  A Lie algebroid over $X$ will then be represented by an $\LA$-groupoid of the form 
\begin{equation*}
\la{A \times \Gamma}{M \times \Gamma}{A}{M},
\end{equation*}
where $\Gamma$ acts on $A$ by Lie algebroid automorphisms.  In other words, the orbifold Lie algebroid cohomology in this case coincides with the equivariant Lie algebroid cohomology $H_\Gamma(A)$.
\end{ex}

\section{Poisson groupoids}\label{sec:cotangent}

Let $(G, \pi) \arrows M$ be a Poisson groupoid.  Mackenzie \cite{mac:dbl} showed that the cotangent prolongation groupoid of $G$ forms an $\LA$-groupoid, as in \eqref{eqn:cotangent}.  We may then apply the $[-1]$ functor to obtain the $Q$-groupoid $[-1]T^*G \arrows [-1]A^*$, where the homological vector field is the Lichnerowicz-Poisson differential $d_\pi$.  The first two rows of the standard cochain complex are of the form
\begin{equation*}
\xymatrix{\vdots & \vdots & \vdots & \\
C^\infty(G) \ar^\delta[u] \ar^{d_\pi}[r] & \vect(G) \ar^{-\delta}[u] \ar^{d_\pi}[r] & \vect^2(G) \ar^\delta[u] \ar^{d_\pi}[r] & \cdots \\
C^\infty(M) \ar^\delta[u] \ar^{d_{A^*}}[r] & \Gamma(A) \ar^{-\delta}[u] \ar^{d_{A^*}}[r] & \bigwedge^2 \Gamma(A) \ar^\delta[u] \ar^{d_{A^*}}[r] & \cdots}
\end{equation*}

We may explicitly describe the first vertical map $\delta^1: \bigwedge^\bullet \Gamma (A) \to \vect^\bullet(G)$.  The face maps $\sigma_i^1$ are given by
\begin{align*}
\left(\sigma_0^1\right)^*(X) &= \linv{X}, \\
\left(\sigma_1^1\right)^*(X) &= \rinv{X},
\end{align*}
for $X \in \bigwedge^\bullet \Gamma(A)$, where $\linv{X}$ and $\rinv{X}$ are, respectively, the left- and right-invariant multivector fields associated to $X$.  Therefore $\ker \delta^1$ consists of multisections $X \in \bigwedge^\bullet \Gamma(A)$ such that $\linv{X} = \rinv{X}$.  It is difficult to describe explicitly the second vertical map, since there does not seem to be a simple description of the cochains that would appear in the higher rows; however, it may be seen \cite{xu:mult} that $\ker \delta^2$ consists of multiplicative multivector fields on $G$.

\begin{ex}
In the case of a Poisson-Lie group $(G, \pi)$, the space of $(p,q)$-cochains is $\bigwedge^p \mathfrak{g} \otimes C^\infty(G^q)$, and $\delta$ is the differential for the smooth group cohomology with coefficients in $\bigwedge \mathfrak{g}$.  If $G$ is compact, the spectral sequence collapses to $\left(\left(\bigwedge \mathfrak{g}\right)^G, d_{\mathfrak{g}^*}\right)$, and we obtain the $G$-invariant Lie algebra cohomology of $\mathfrak{g}^*$.
\end{ex}

\section{Vacant $\LA$-groupoids}\label{sec:vacant}

An $\LA$-groupoid (\ref{eq:la2}) is said to be \emph{vacant} if it has trivial core\footnote{See \cite{mac:dblie1} for a description of the core of an $\LA$-groupoid and its induced Lie algebroid structure.} or, equivalently, if the induced map $\Omega \to G \times_M A$ is a diffeomorphism.  Since the fibres of this map are vector spaces, it suffices to count dimensions in order to check if an $\LA$-groupoid is vacant.  Examples of vacant $\LA$-groupoids include $\LA$-groupoids representing Lie algebroid structures over orbifolds (\S\ref{sec:orb}) and the cotangent groupoid (\S\ref{sec:cotangent}) of a Poisson-Lie group.  

Mackenzie \cite{mac:dblie1} has shown that every vacant $\LA$-groupoid is isomorphic to a ``matched pair'' $\LA$-groupoid $G \bowtie A$, constructed out of a compatible pair of actions of $G$ and $A$ on each other.  We will review the matched pair construction and then describe the double complexes arising from vacant $\LA$-groupoids.

First recall the notions of groupoid and Lie algebroid actions \cite{mac:lie}.  If $G \arrows M$ is a Lie groupoid, then a left action of $G$ on a vector bundle\footnote{One can, of course, define Lie groupoid actions on general fibre bundles, but vector bundle actions are sufficient for the present purposes.} $E \to M$ is a linear map $s^*(E) \to E$, $(g,v) \mapsto g(v)$, such that the diagram
\begin{equation*}
\xymatrix{s^*(E) \ar[r] \ar[d] & G \ar^t[d]\\ E \ar[r] & M}
\end{equation*}
commutes and such that $g(h(v))=(gh)(v)$ for all $(g,h,v) \in (s \circ m)^*(E)$.

A left action of $G$ on $E$ induces a Lie groupoid structure on $s^*(E) \arrows E$, where the source and target maps are given by
\begin{align*}
\tilde{s}(g,v) &= v,\\
\tilde{t}(g,v) &= g(v).
\end{align*}
A pair $\left( (g,v),(h,w) \right)$ is composable if $v = h(w)$.  It follows, in particular, that $g$ and $h$ are composable elements of $G$, and the multiplication is then defined as $(g,v) \cdot (h,w) = (gh,w)$.  We leave descriptions of the identity and inverse maps as an exercise for the reader.

Now let $A \to M$ be a Lie algebroid, and let $P \stackrel{\pi}{\to} M$ be a submersion.  Then a (right) action of $A$ on $P$ is a Lie algebra homomorphism $\tilde{\rho}: \Gamma(A) \to \vect(P)$ that lifts the anchor map $\rho: \Gamma(A) \to \vect(M)$, in the sense that $\tilde{\rho}(X)$ is $\pi$-related to $\rho(X)$ for all $X \in \Gamma(A)$.  

An action of $A$ on $P$ induces a Lie algebroid structure on $\pi^*(A) \to P$, as follows.  The action map $\tilde{\rho}$ may be extended by $C^\infty(P)$-linearity to obtain the map (which we will also denote $\tilde{\rho}$)
\begin{equation*}
\tilde{\rho}: \Gamma(\pi^*(A)) = C^\infty(P) \otimes \Gamma(A) \to \vect(G),
\end{equation*}
which is the anchor map for the induced Lie algebroid.  The Lie bracket of sections is defined by setting
\begin{equation*}
[1 \otimes X, 1 \otimes Y] = 1 \otimes [X,Y]
\end{equation*}
for $X,Y \in \Gamma(A)$ and extending by the Leibniz rule.  The Jacobi identity follows from the Jacobi identity for the bracket on $\Gamma(A)$ and the fact that $\tilde{\rho}$ is a Lie algebra homomorphism.

Let $G \arrows M$ be a Lie groupoid and let $A \to M$ be a Lie algebroid (not necessarily the Lie algebroid of $G$), such that $G$ is equipped with an action on $A$, and $A$ is equipped with an action on $G \stackrel{s}{\to} M$.  Then, based on the above discussion, we may form a square
\begin{equation}\label{vacantla}
\la{s^*(A)}{G}{A}{M},
\end{equation}
where the horizontal and vertical sides are respectively Lie algebroids and Lie groupoids.  It is automatically true that $\tilde{s}$ is a Lie algebroid homomorphism; if $\tilde{t}$ and $\tilde{m}$ are also Lie algebroid homomorphisms, then (\ref{vacantla}) is an $\LA$-groupoid.

\begin{dfn}A \emph{groupoid-algebroid matched pair} is a Lie groupoid $G \arrows M$ and a Lie algebroid $A \to M$ equipped with mutual actions such that $\tilde{t}$ and $\tilde{m}$ are Lie algebroid morphisms.
\end{dfn}

\begin{rmk}The reader may wish to see Mackenzie \cite{mac:dblie1} for a more concrete set of compatibility conditions for groupoid-algebroid matched pairs.\end{rmk}

Suppose that $(G,A)$ is a matched pair.  We apply the $[-1]$ functor to (\ref{vacantla}) to obtain the $Q$-groupoid $s^*([-1]A) \arrows [-1]A$, whose homological vector field we denote as $d_{\tilde{A}}$.  The algebra of functions on $s^*([-1]A)$ is $C^\infty(G) \otimes_s \bigwedge \Gamma(A^*)$.  For the higher groupoid cochains, we observe that
\begin{equation*}\left(s^*([-1]A)\right)^{(q)} = (s \circ p^q_q)^*([-1]A),
\end{equation*}
where $p^q_q: G^{(q)} \to G$ is the projection map onto the last component.  So
\begin{equation*}C^\infty\left(\left(s^*([-1]A)\right)^{(q)}\right) = C^\infty(G^{(q)}) \otimes_{s \circ p^q_q} \bigwedge \Gamma(A^*),
\end{equation*}
and the space of $(p,q)$-cochains for the double complex is therefore 
\begin{equation}\label{vacantchains}
C^{p,q}(s^*([-1]A)) = C^\infty(G^{(q)}) \otimes_{s \circ p^q_q} \bigwedge^p \Gamma(A^*).
\end{equation}

As is already clear from (\ref{vacantchains}), the double complex intertwines the Lie groupoid cohomology of $G$ and the Lie algebroid cohomology of $A$.  If $G$ has compact $t$-fibres, then the spectral sequence collapses at the first stage to the $G$-invariant Lie algebroid complex of $A$.

\begin{ex}[\cite{mac:reduction,mythesis}]\label{ex:lu}
Let $M$ be a Poisson manifold, and let $H$ be a Poisson-Lie group with a Poisson action $s: M \times H \to M$.  Then there an action of the groupoid $M \times H$ on $T^*M$, given by the map $\tilde{t}: s^*(T^*M) \to T^*M$, sending $(x,g,\eta)$, where $\eta \in T^*_{xg}M$, to $r_g^*\eta$, where $r_g$ is the map given by right-multiplication by $g$.  Additionally, there is an action of $T^*M$ on $M \times H$, given by the map $\tilde{\rho}: \Omega^1(M) \to \vect(M \times G)$, $\alpha \mapsto \tilde{\pi}^\sharp(s^*\alpha)$, where $\tilde{\pi}^\sharp: \Omega^1(M \times H) \to \vect(M \times H)$ arises from the product Poisson structure on $M \times H$.  This pair of actions gives $(M \times H, T^*M)$ the structure of a matched pair, so that
\begin{equation}\label{eqn:lu}
\la{s^*(T^*M)}{M \times H}{T^*M}{M}
\end{equation}
is an $\LA$-groupoid.  

This $\LA$-groupoid was introduced by Mackenzie \cite{mac:reduction} in order to describe a general procedure for Poisson reduction.  By applying the Lie functor to the vertical sides of (\ref{eqn:lu}), one obtains the double Lie algebroid corresponding to the matched pair Lie algebroid structure, due to Lu \cite{lu}, on $(M \times \mathfrak{h}) \oplus T^*M$.

If $H$ is compact, then the $\LA$-groupoid cohomology of (\ref{eqn:lu}) is equal to the $H$-invariant Poisson cohomology of $M$.
\end{ex}

\bibliographystyle{abbrv}
\bibliography{bibio}

\end{document}